\documentclass[12pt,a4paper]{article}
\usepackage{hyperref}
\hypersetup{pdfstartview=FitH}
\hypersetup{pdfpagemode=FitWidth}
\usepackage{amssymb}
\usepackage{amsmath}
\usepackage{amsthm}
\usepackage{times} 
\newcommand{\n}{\newcommand} 
\n{\oo}{\operatorname}
\n{\am}{\begin{pmatrix}}\n{\zm}{\end{pmatrix}}
\n{\mb}{\mathbb}\n{\mc}{\mathcal}
\n{\bb}{\bigskip}
\n{\G}{\mc{G}}
\n{\F}{\mc{F}}
\n{\g}{\gamma}\n{\Ga}{\Gamma}\n{\HH}{\mc{H}}
\n{\I}{\mb{I}}\n{\II}{\mc{I}}
\n{\Q}{\mb{Q}}
\n{\Z}{\mb{Z}}
\newcommand{\ii}{\hskip 1em\relax}
\newcommand{\q}{\quad} 

\begin{document}
\begin{center}
{\Huge Hurwitz's Freeness Property}
\end{center}\bb\bb 

In \cite{g} Gauss devised an algorithm to solve in integers the equation 

\begin{equation}\label{m}
a\,x^2+2\,b\,x\,y+c\,y^2=m,
\end{equation}

where $a,b,c,m$ are given integers. \bb 

\ii Consider the groups

$$G:=\oo{G L}(2,\Z)/\{\pm1\}
\supset\oo{S L}(2,\Z)/\{\pm1\}=:H.$$

Let $\Delta$ be an integer, let $F$ be the set all quadratic forms 

\begin{equation}\label{f}
f(X,Y):=a\,X^2+2\,b\,X\,Y+c\,Y^2
\end{equation}

with $a,b,c$ integers and $b^2-a c=\Delta$, and let $\F$ be the groupoid attached to the natural action of $H$ on $F$. In Article 169 of \cite{g} Gauss reduces the solution of (\ref{m}) to the computation of certain hom-sets in $\F$ (see below).\bb 

\ii Assume from now on that $\Delta$ is a fixed positive nonsquare integer.\bb 

\ii The group $G$ acts on the set $\I$ of irrational real numbers by linear fractional transformations. Let $\II$ be the corresponding groupoid, and $\I_\Delta$ the set of real numbers 

$$\frac{-b-\sqrt{\Delta}}{a}$$ 

where $f$ as in (\ref{f}) runs over the elements of $F$. Then $H$ preserves $\I_\Delta$, and we can form the restricted groupoid $\II_\Delta^H$ issued from the action of $H$ on $\I_\Delta$.\bb 

\ii In Section 73 of \cite{d} Dirichlet notes that above formula gives a canonical groupoid isomorphism from $\F$ to $\II_\Delta^H$.\bb 

\ii In Section 63 of \cite{h} Hurwitz shows that the groupoid $\II$ is free over one of its sub oriented graph, giving a very simple description of the hom-sets of $\II_\Delta^H$, which Dirichlet had identified to the hom-sets of $\F$, whose computation Gauss had reduced the solution of (\ref{m}) to.\bb  

\ii We wish to phrase Hurwitz's statement in today's language.\bb 

\ii Say that the {\bf derivative} $x'$ of a point $x$ in $\I$ is the inverse of its fractional part, let $g_x$ be the image in $G$ of 

$$\am\lfloor x\rfloor&1\\1&0\zm,$$

so that we have $x=g_x\,x'$, let $\g(x)$ be the corresponding morphism in $\II$ from $x'$ to $x$, and let $\Ga$ be the sub oriented graph of $\II$ whose vertices are the points of $\I$ and whose arrows are the $\g(x)$.\bb 

\ii Then $\II$ is the groupoid freely generated by $\Ga$ in the following sense.\bb 

\ii Let $\varphi$ be an oriented graph morphism from $\Ga$ into any groupoid $\G$. Then $\varphi$ extends uniquely to a groupoid morphism from $\II$ to $\G$.\bb 

\ii The groupoid $\II$ has a very simple structure, which can be described as follows. To ease notation put $x_i:=x^{(i)}$.\bb 

\ii Let $g$ be a nontrivial morphism in $\II$ from $x$ to $y$. Then there is a unique pair $(i,j)$ of nonnegative integers satisfying $x_i=y_j$, 

$$g=\g(y_0)\,\cdots\,\g(y_{j-1})\ \g(x_{i-1})^{-1}
\,\cdots\,\g(x_0)^{-1},$$

and $x_{i-1}\not=y_{j-1}$ if $i$ and $j$ are positive.\bb 

\ii The composition of two such elements is tedious but easy to compute.\bb 

\ii Let $x$ be in $\I$. Recall that the sequence $(x_i)$ is eventually periodic if and only if $x$ has degree 2 over $\Q$. This makes $\II_\Delta^H$ computable. In particular the stabilizer in $H$ of $f$ in $F$ is infinite cyclic. However, $\II$ is highly uncomputable.\bb 

\centerline{*}\bb

\ii We said above that Gauss reduced the solution of (\ref{m}) to the computation of certain hom-sets in $\F$. Let's be more precise. (We only indicate some of the main statements, directing the reader to \cite{h} for a full treatment.) Assume that $m$ is nonzero.\bb 

\ii Say that a solution of (\ref{m}) is a {\bf representation of $m$ by} $f$ ($f$ being given by (\ref{f})), and that such a representation is {\bf proper} of $X$ and $Y$ are relatively prime. It clearly suffices to describe the set $P$ of proper representations of $m$ by $f$.\bb

\ii As a general notation, write $[a,b,c]$ for the form (\ref{f}). Let $n$ be in $N$. Put

$$\ell_n:=\frac{n^2-\Delta}{m}\q,\q f_n:=[m,n,\ell_n],$$ 

form the set $S_n$ of those substitutions $h$ in $\oo{S L}(2,\Z)$ which satisfy $f h=f_n$, and let $u_n$ be the map from $S_n$ to $P$ attaching to $h\in S_n$ its first column.\bb

\ii Then the $u_n$ induce a bijection form the disjoint union of the $S_n$ onto $P$.

\vfill Pierre-Yves Gaillard\hfill\tiny hurwitz.freeness.080924, Wed Sep 24 07:54:35 CEST 2008. 

\end{document}